\documentclass{amsart}

\DeclareMathSymbol{\twoheadrightarrow}  {\mathrel}{AMSa}{"10}

\def\Q{{\mathbf Q}}
\def\Z{{\mathbf Z}}

\def\F{{\mathbf F}}

\def\W{{V}}

\def\Gal{\mathrm{Gal}}

\def\Aut{\mathrm{Aut}}
\def\Hom{\mathrm{Hom}}
\def\I{{\mathcal I}}
\def\II{{{\mathcal I}_{v,X}}}

\def\r{{\mathfrak G}}

\def\GL{\mathrm{GL}}

\def\Sp{\mathrm{Sp}}
\def\M{\mathcal{M}}
\def\N{\mathcal{N}}
\def\L{\mathcal{L}}

\def\dim{\mathrm{dim}}
\def\rk{\mathrm{rk}}

\def\T{{\mathbf T}}
\def\B{{\mathbf B}}

\def\invlim{{\displaystyle{\lim_{\leftarrow}}}}

\newtheorem{thm}{Theorem}[section]
\newtheorem{lem}[thm]{Lemma}
\newtheorem{cor}[thm]{Corollary}
\newtheorem{prop}[thm]{Proposition}
\theoremstyle{definition}
\newtheorem{defn}[thm]{Definition}

\newtheorem{rem}[thm]{Remark}

\title[Subgroups of inertia groups]
{Subgroups of inertia groups arising from abelian varieties}

\author[A.\ Silverberg]{A.\ Silverberg}
\address{Department of Mathematics, Ohio State University, 
231 W.\ 18 Avenue,
Columbus, Ohio 43210--1174, USA}
\email{silver\char`\@math.ohio-state.edu}

\author[Yu. G. Zarhin]{Yu. G. Zarhin}
\address{Department of Mathematics, Pennsylvania State University, 
University Park, PA 16802, USA}
\email{zarhin\char`\@math.psu.edu}

\begin{document}

\begin{abstract}
Given an abelian variety over a field with a discrete valuation,
Grothendieck defined a certain open normal subgroup of the 
absolute inertia group. This  
subgroup encodes information on the extensions over which the 
abelian variety acquires semistable reduction. 
We study this subgroup, and use it to obtain information on the
extensions over which the abelian variety acquires
semistable reduction.
\end{abstract}

\maketitle
\baselineskip=15pt

\section{Introduction}

Suppose $X$ is an abelian variety over a field $F$, and
$v$ is a discrete valuation on $F$. Fix an extension ${\bar v}$
of $v$ to a separable closure $F^{s}$ of $F$, and write
$\I_{v}$ for the inertia subgroup in $\Gal(F^{s}/F)$ for ${\bar v}$.
In \cite{SGA} (see pp.~354--355), Grothendieck defined a subgroup $\I'$ of 
$\I_{v}$ with the property that $X$ has semistable reduction at 
the restriction $w$ of ${\bar v}$ to a finite separable extension of 
$F$ if and only if $\I_{w} \subseteq \I'$. 
In particular, if $F_{v}^{nr}$ denotes the maximal unramified
extension of the completion of $F$ at $v$, then $\I'$ cuts out
the smallest Galois extension of $F_{v}^{nr}$ over which $X$ has
semistable reduction. 
We denote the group $\I'$ by $\II$ because of its dependence
on $X$ and $v$. 

In \S\ref{properties} we give some 
properties of the group $\II$. 
We show that the Zariski closure of its image under the $\ell$-adic
representation (for $\ell$ different from the residue characteristic)
coincides with the identity connected component of
the Zariski closure of the image of $\I_{v}$. 
The proofs of the results in \S\ref{properties}
are in the spirit of \cite{Compositio},
where we dealt with connectedness questions for Zariski closures
of images of $\ell$-adic representations. 

In \S\ref{bounds} we show that the finite group 
$G_{v,X} = \I_{v}/\II$ injects into 
$\GL_{t-t_{v}}(\Z) \times \Sp_{2(a-a_{v})}(\Z_{\ell})$
for all but finitely many primes $\ell$, 
where $t_{v}$ and $a_{v}$ (respectively, $t$ and $a$)
are the toric and abelian ranks of the special fiber of
the N\'eron model of $X$ at $v$ (respectively, at an
extension of $v$ over which $X$ has semistable reduction).
Here, the projection onto the first factor is independent of $\ell$,
and the characteristic polynomial of the projection onto the
second factor has integer coefficients independent of $\ell$.
The group $G_{v,X}$ was introduced by Serre in the case
of elliptic curves in 
\S 5.6 of \cite{Serre72} (where it was called $\Phi_{p}$).

In \S\ref{order} we obtain divisibility bounds on the order of $G_{v,X}$, 
and in \S\ref{applic} we deduce results on 
semistable reduction of abelian varieties.
Bounds on the prime divisors and the exponent
of $\#G_{v,X}$ were obtained
by Lorenzini (see Proposition 3.1 of \cite{Lorenzini}).
In particular, the bound on $Q_{v,X}$ in Corollary \ref{hdcor}
was essentially obtained by Lorenzini.

The paper continues our earlier work on semistable reduction
of abelian varieties (see \cite{semistab}
and \cite{smalldeg}). The proofs are heavily influenced by 
the fundamental results of Grothendieck and Serre.

The authors would like to thank Karl Rubin for helpful
conversations, and NSA and NSF for financial support.

\section{Notation and preliminaries}
\label{prelim}

If $K$ is a field, write $K^s$ for a separable closure.
If $G$ is an algebraic group, let $G^{0}$
denote its identity connected component. 
Let $\varphi$ denote the Euler $\varphi$-function, let
$\zeta_{M}$ denote a primitive $M$-th root of unity,
and let $\F_{\ell}$ denote the finite field with $\ell$ elements.

If $Y$ is a commutative algebraic group over a field $K$
(e.g., an abelian variety or an algebraic torus), let 
$Y_{n}$ denote the kernel of multiplication by $n$ in $Y(K^{s})$, 
let 
$T_{\ell}(Y) = \invlim Y_{\ell^{n}}$, and let $V_{\ell}(Y) = 
T_{\ell}(Y) \otimes_{\Z_{\ell}} \Q_{\ell}$. 
For example, if ${\mathbf G}_{m}$ is the multiplicative group,
then $\Z_\ell(1):=T_{\ell}({\mathbf G}_{m})$ is a free 
$\Z_{\ell}$-module of rank $1$. 

Throughout this paper, 
$X$ is a $d$-dimensional abelian variety over a field $F$, 
$v$ is a discrete valuation on $F$ of residue characteristic $p \ge 0$,
and $\ell$ is a prime not equal to $p$.
Fix an extension ${\bar v}$ of $v$ to 
$F^s$. 
If $w$ is the restriction of ${\bar v}$ to a finite separable 
extension $L$ of $F$, let $\I_{w}$ denote the inertia subgroup in 
$\Gal(F^{s}/L)$ for ${\bar v}$,
let $X_{w}$ denote the special fiber of the N\'eron model
of $X$ at $w$, and
let $\T_{w}$ denote the maximal subtorus of $X_{w}$.
Let
$$\rho_{\ell}:\Gal(F^s/F)\to\GL(T_\ell(X))$$
denote the $\ell$-adic representation.
Let 
$\r$ denote the Zariski closure of $\rho_{\ell}(\I_{v})$ in
$\GL(V_{\ell}(X))$. 
We will make repeated use of the following result.

\begin{thm}[Galois Criterion for Semistable Reduction]
\label{galcrit}
$X$ has semistable reduction at $v$ if and only if 
$\I_{v}$ acts unipotently on $V_\ell(X)$.
\end{thm}

\begin{proof}
See Proposition 3.5 and Corollaire 3.8 of \cite{SGA} 
and Theorem 6 on p.~184 of \cite{BLR}.
\end{proof}

\begin{thm}
\label{conn}
The Zariski closure $\r$ of $\rho_{\ell}(\I_v)$ is connected
if and only if $X$ has semistable reduction at $v$.
\end{thm}

\begin{proof}
See Theorem 5.2 of \cite{Compositio}; see also Remarque~1 on
p.~396 of \cite{Motives}. 
\end{proof}

Suppose $\lambda$ is a polarization on $X$.
Then $\lambda$ gives rise to a non-degenerate, alternating,
$\Gal(F^{s}/F)$-equi\-var\-i\-ant, $\Z_\ell(1)$-valued 
pairing on $T_{\ell}(X)$ (see \S\S 1.0 and 2.5 of \cite{SGA}).
Since $\I_v$ acts trivially on $\Z_\ell(1)$, we obtain
a non-degenerate, alternating, $\I_v$-invariant pairing
$$E_\lambda : T_\ell(X) \times T_\ell(X) \to \Z_\ell.$$
Then $E_\lambda$ is perfect if and only if $\deg(\lambda)$
is not divisible by $\ell$ (see (2.5.1) and \S 1.0 of
\cite{SGA}).
Let $\perp$ denote the orthogonal complement with respect to 
$E_\lambda$.

\begin{prop}
\label{serretateeqn}
\begin{enumerate}
\item[(i)]
We can identify
$$T_\ell(X_v)=T_\ell(X)^{\I_v} \quad \text{ and } \quad
V_\ell(X_v)=V_\ell(X)^{\I_v}.$$
\item[(ii)] $T_\ell(\T_v)=T_\ell(X)^{\I_v} \cap (T_\ell(X)^{\I_v})^\perp
= T_\ell(X_v) \cap T_\ell(X_v)^{\perp}$.
\item[(iii)] $T_\ell(X_v)/T_\ell(\T_v)$ is a free $\Z_\ell$-module.
\end{enumerate}
\end{prop}

\begin{proof}
For (i), see Lemma 2 on p.~495 of \cite{SerreTate}.
For (ii), use (i) and Grothendieck's Orthogonality Theorem 
(see (2.5.2) of \cite{SGA}).
For (iii), see (2.1.6) of \cite{SGA}.
\end{proof}

\begin{prop}
\label{orthogeqn}
Suppose $L$ is a finite separable extension of $F$,  $w$
is the restriction of ${\bar v}$ to $L$, and
$X$ has semistable reduction at $w$.
Then$$ T_\ell(\T_w) = T_\ell(X_w)^\perp
= (T_\ell(X)^{\I_w})^\perp
 \subseteq  
T_\ell(X)^{\I_w} = T_\ell(X_w).$$
\end{prop}

\begin{proof}
See Proposition 3.5 of \cite{SGA}.
\end{proof}

\section{Linear Algebra}

We will use the following linear algebra facts 
in Theorem \ref{hd} below.

\begin{defn}
Suppose $R$ is a principal ideal domain, and $\M_1$ 
and $\M_2$ are free $R$-modules.
A bilinear form $e : \M_1 \times \M_2 \to R$ 
is called a {\it perfect} pairing if the natural 
homomorphisms 
$\M_1 \to \Hom(\M_2,R)$ and $\M_2 \to \Hom(\M_1,R)$
are bijective.
If $\L$ is a submodule of $\M_1$ (resp., $\M_2$), we write 
$\L^{\perp}$ for the orthogonal complement of $\L$
with respect to $e$ in $\M_2$ (resp., $\M_1$).
\end{defn}

\begin{rem}
\label{perf}
Suppose $R$ is a principal ideal domain,  
$\M$  is a free $R$-module of rank $2n$, and
$e : \M \times \M \to R$ is an alternating bilinear form.
If $e$ is perfect, then $\Aut(\M,e) \cong \Sp_{2n}(R)$,
where  
$\Sp_{2n}(R)$ denotes the group of $2n \times 2n$ symplectic
matrices over $R$ (see \S{5} of \cite{Bourbaki}).
\end{rem}

\begin{prop}
\label{linalg}
Suppose $\ell$ is a prime number, $G$ is a finite group whose order
is not divisible by $\ell$, $V$ is a finite-dimensional
$\Q_\ell$-vector space with a linear $G$-action, and
$e:V \times V \to \Q_\ell$ is a $G$-invariant
non-degenerate alternating (resp., symmetric) bilinear form. 
Suppose $\M$ is a $G$-stable $\Z_\ell$-lattice 
in $V$ (there always exist such). 
Then there exists a perfect $G$-invariant alternating 
(resp., symmetric) bilinear
form $e':\M \times \M \to \Z_\ell$.
\end{prop}

\begin{proof}
First note that the result is true when $e$ is symmetric
and $\dim(V)=1$ and when $e$ is alternating and
$\dim(V)=2$, by scaling $e$.

We may assume that
$e(\M,\M)=\Z_{\ell}$ (by replacing 
$e$ by $\ell^i e$ for a suitable $i$, 
if necessary).  
Then $e$ induces a non-zero alternating (resp., symmetric) 
bilinear form $\bar{e}$ on $\M/\ell \M$. 
By Nakayama's Lemma,
if $\bar{e}$ is non-degenerate then $e$ is perfect
and we are done.
Assume $\bar{e}$ is degenerate, and let
$\bar{\M_1}=\ker(\bar{e})$. Then $\bar{\M_1}$ is a proper
subset of $\M/\ell \M$, 
since $e(\M,\M)=\Z_{\ell}$.
Since $\#G$ is not divisible by $\ell$,
there is a $G$-invariant splitting
        $\M/\ell \M= \bar{\M_1} \oplus \bar{\M_2}$,
where $\bar{\M_2}$ is a non-zero $G$-invariant subspace of $\M/\ell \M$,
which can be lifted to a $G$-invariant splitting of 
$\Z_{\ell}$-lattices 
        $\M=\M_1 \oplus \M_2$
with $\M_1/\ell \M_1=\bar{\M}_1$ and $\M_2/\ell \M_2=\bar{\M}_2$
(see \S 15.5 and Corollary 1 of \S 14.4  of \cite{serrereps}).
Denote by $e_2$ the restriction of $e$ to
$V_2:=\M_2\otimes_{\Z_\ell} \Q_{\ell}$, and 
let $V_1$ denote the orthogonal complement of $V_2$ in $V$ with respect
to $e$. The restriction of $e$ to $V_1$ is non-degenerate, and
the restriction of $e$ to $\M_2$ is perfect.
We obtain a $G$-invariant orthogonal splitting
$V= V_1 \oplus V_2$. 
Replace $\M_1$ by its (isomorphic) image in $V_1$ under
the projection map from $V_1 \oplus V_2$ to $V_1$. 
Since $\dim(V_1) < \dim(V)$, we obtain inductively a
perfect
$G$-invariant alternating (resp., symmetric) form  $e_1$ on $\M_1$.
Let $\M=\M_1 \oplus \M_2$ and $e'=e_1 \oplus e_2$.
\end{proof}

\begin{lem}
\label{perp}
Suppose $R$ is a principal ideal domain,  $\M_1$ 
and $\M_2$ are free $R$-modules, $\L$ is a submodule
of $\M_1$, $\M_1/\L$ is torsion-free, and
$e : \M_1 \times \M_2 \to R$ 
is a {\it perfect} pairing.
Then:
\begin{enumerate}
\item[(i)] $(\L^{\perp})^{\perp} = \L$,
\item[(ii)] the natural map 
$\Hom(\M_1,R) \to \Hom(\L,R)$ is surjective,
\item[(iii)] the induced form
$\L \times \M_2/\L^{\perp} \to R$ is a perfect pairing.
\end{enumerate} 
\end{lem}

\begin{proof}
Clearly, $\L \subseteq (\L^{\perp})^{\perp}$.  
Let $K$ denote the fraction field of $R$.
By dimension arguments we have
$\L \otimes K = ((\L \otimes K)^{\perp})^{\perp} =
(\L^{\perp})^{\perp}  \otimes K$.
Therefore, $(\L^{\perp})^{\perp}/\L$ is torsion.
Since $\M_1/\L$ is torsion-free, we obtain (i).
Further, since $\M_1/\L$ is torsion-free, 
$\L$ is a direct summand of $\M_1$, and therefore we have (ii).

Since $e$ is perfect, we can identify $\Hom(\M_2,R)$ with $\M_1$.
Under this identification,
$\Hom(\M_2/\L^{\perp},R) = (\L^{\perp})^{\perp} = \L$,
by (i).  
Since $e$ is perfect we can identify $\M_2$ with $\Hom(\M_1,R)$.
The natural injection 
$\M_2/\L^{\perp} \hookrightarrow \Hom(\L,R)$ is surjective 
since $\M_2 = \Hom(\M_1,R) \to \Hom(\L,R)$ 
is surjective by (ii).
\end{proof}

\begin{prop}
\label{linalg1}
Suppose $R$ is a principal ideal domain, $\M$ is a free $R$-module
of finite rank,  $e : \M \times \M \to R$ is an alternating
(resp., symmetric) perfect
pairing, $\N$ is a submodule of $\M$, and $\M/\N$ is
torsion-free. 
Assume that $\N^{\perp} \subseteq \N$. 
Then $${\tilde e}(a+\N^{\perp},b+\N^{\perp})=e(a,b)$$
defines an alternating (resp., symmetric) perfect  
pairing
${\tilde e} : \N/\N^{\perp} \times \N/\N^{\perp} \to R$.
\end{prop}

\begin{proof}
Applying Lemma \ref{perp}iii with $\M_1 = \M_2 = \M$ and $\L = \N$ gives
an alternating (resp., symmetric) perfect pairing 
$e' : \N \times \M/\N^{\perp} \to R$.  
Now applying Lemma \ref{perp}iii to $e'$ with 
$\M_1 = \M/\N^{\perp}$, $\M_2 = \N$, and $\L = \N/\N^{\perp}$,
we obtain the desired result.
\end{proof}

\section{Properties of $\II$}
\label{properties}

We will give a different definition for $\II$  
than Grothendieck did, and will then show that the two definitions
are equivalent. Grothendieck's definition coincides with (ii) of
Theorem \ref{main} below. 

\begin{defn}
Define 
$\II$ to be the kernel of the natural surjective homomorphism
$\I_{v} \to \r/\r^{0}$. 
Define  
$G_{v,X}  = \I_{v}/\II$.
\end{defn}

\begin{thm}
\label{main}
$\II$ is an open normal subgroup of $\I_{v}$ 
which enjoys the following properties:
\begin{enumerate}
\item[(i)] $\II$ is the largest open
subgroup of $\I_{v}$ such that the Zariski closure of its
image under $\rho_{\ell}$ is $\r^{0}$.
\item[(ii)] $\II = 
\{\sigma \in \I_{v} : 
\sigma \text{ acts unipotently on } V_{\ell}(X) \}$.
\item[(iii)] If $L$ is a finite separable extension of $F$
and $w$ is the restriction of ${\bar v}$ to $L$, then the
following are equivalent:
\begin{enumerate}
\item[(a)] $X$ has semistable reduction at $w$,
\item[(b)] $\I_{w} \subseteq \II$,
\item[(c)] 
the Zariski closure of $\rho_{\ell}(\I_{w})$ is $\r^{0}$.
\end{enumerate} 
\item[(iv)] $\II = \I_{v}$ if and only if $X$ has semistable
reduction at $v$.
\item[(v)] $\II$ is independent of the choice of $\ell$.
\end{enumerate} 
\end{thm}

\begin{proof}
By the definition of $\II$, 
it is an open normal subgroup of $\I_{v}$,
and it is the largest open
subgroup of $\I_{v}$ such that the Zariski closure of its
image under $\rho_{\ell}$ is $\r^{0}$. 

Suppose $L$ is a finite separable extension of $F$, and
$w$ is the restriction of ${\bar v}$ to $L$.
Let $\r_{w}$ denote
the Zariski closure of $\rho_{\ell}(\I_{w})$,
and let $\W = V_{\ell}(X)^{\I_{w}}$.

Suppose first that
$X$ has semistable reduction at $w$.
Then $\r_{w}$ is connected, by Theorem \ref{conn}. 
Therefore, $\I_{w} \subseteq \II$ and
$\r_{w} = \r^{0}$.
Since $\r^{0}$ ($ = \r_{w}$) is a
normal subgroup of $\r$, it follows that $\W$ is stable under
$\I_{v}$.  
Next we will show that 
$$\II = 
\{\sigma \in \I_{v} : 
\sigma \text{ acts unipotently on } V_{\ell}(X) \}$$
$$= 
\{\sigma \in \I_{v} : 
\sigma \text{ acts unipotently on } \W \}.$$
Since $\r^{0} = \r_{w}$,
$\r^{0}$ acts as the identity on $\W$ and on $V_{\ell}(X)/\W$.
Therefore, every element of $\II$ acts unipotently on $V_{\ell}(X)$,
and therefore on $\W$.
To show the reverse inclusions, 
suppose $g \in \I_{v}$ and $g$ acts unipotently on $\W$.
By Proposition \ref{orthogeqn} (and after tensoring with
$\Q_\ell$), 
$\W^{\perp} \subseteq  \W$.
Since $\W^{\perp}$ is the dual of
$V_{\ell}(X)/\W$, it follows that $g$ acts unipotently on
$V_{\ell}(X)/\W$, and therefore acts unipotently on
$V_{\ell}(X)$. By Proposition 2.5 of \cite{Compositio},
$g \in \II$. We therefore obtain the desired equalities.

By (i), if
$\r^{0} = \r_{w}$, then $\I_{w} \subseteq \II$.
By (ii) and Theorem \ref{galcrit}, 
if $\I_{w} \subseteq \II$ then $X$ has semistable
reduction at $w$. We therefore have (iii). 
We easily deduce (iv) from (iii).

By Th\'eor\`eme 4.3 of \cite{SGA},
if $\sigma \in \I_{v}$ then the characteristic polynomial
of $\rho_{\ell}(\sigma)$ is independent of $\ell$.
By (ii), 
$$\II =
\{\sigma \in \I_{v} : \text{ the characteristic polynomial
of } \rho_{\ell}(\sigma) \text{ is } (x-1)^{2d}\}.$$
Therefore, $\II$ is independent of $\ell$.
\end{proof}

\begin{prop}
\label{Wunipot}
If $L$ is a finite separable extension of $F$, and $X$ has
semistable reduction at 
the restriction $w$ of ${\bar v}$ to $L$, then 
\begin{enumerate} 
\item[(i)] $\II = 
\{\sigma \in \I_{v} : 
\sigma \text{ acts unipotently on } V_{\ell}(X)^{\I_{w}} \}$,
\item[(ii)]  $G_{v,X}$ acts faithfully on 
$T_{\ell}(X)^{\I_{w}}/T_{\ell}(X)^{\I_{v}}$,
\item[(iii)]  $T_{\ell}(X)^{\I_{w}} = T_{\ell}(X)^{\II}$.
\end{enumerate} 
\end{prop}

\begin{proof}
The proof of Theorem \ref{main} included a proof of (i), and
easily implies (ii). For (iii), let 
$\r_{w}$ denote
the Zariski closure of $\rho_{\ell}(\I_{w})$, and
note that 
$T_{\ell}(X)^{\II} = T_{\ell}(X)^{\r^{0}} = 
T_{\ell}(X)^{\r_{w}} = T_{\ell}(X)^{\I_{w}}$,
by Theorem \ref{main}.
\end{proof}

\section{Restrictions on $G_{v,X}$}
\label{bounds}

\begin{rem}
\label{cyclicrem}
Note that 
$G_{v,X}$ is isomorphic to the group of connected components of $\r$.
If $p=0$ then $G_{v,X}$ is cyclic, and if
$p>0$ then $G_{v,X}$ is an extension of a cyclic group of order
prime to $p$ by a $p$-group 
(as can be seen by replacing $F$ by the maximal unramified
extension of the completion of $F$ at $v$, looking at the
extension cut out by $\II$, taking its maximal tamely ramified
subextension, and applying \S8 of \cite{Frohlich}).
In particular, $G_{v,X}$ is solvable, and has a unique 
Sylow-$p$-subgroup if $p>0$. 
\end{rem}

The group $G_{v,X}$ does not change if we replace $F$ by an
extension unramified at $v$. Therefore, in this section we
can and do replace $F$ by the maximal unramified
extension of the completion of $F$ at $v$. Then $\II$ cuts out
the smallest Galois extension $L$ of $F$ over which $X$ has
semistable reduction. Let $w$ denote the restriction of ${\bar v}$
to $L$. 
Then $\II = \I_{w}$. 
Therefore, $\I_w$ is an open normal subgroup of $\I_v$ of finite index, 
and  
$T_\ell(X_w)=T_\ell(X)^{\I_w}$ is $\I_v$-stable. 
By Proposition \ref{serretateeqn}i we have 
$$T_{\ell}(X_{v}) = T_{\ell}(X)^{\I_{v}} \subseteq 
T_{\ell}(X)^{\I_{w}} = T_{\ell}(X_{w})$$
as $\I_v$-modules.

Over an algebraic closure of the residue field, there are 
exact sequences
$$0 \to \T_{w} \to X_{w}^{0} \to \B_{w} \to 0, \qquad 
0 \to U_{v} \times \T_{v} \to X_{v}^{0} \to \B_{v} \to 0,$$
where $\B_{w}$ and $\B_{v}$ are abelian varieties  
and $U_{v}$ is a unipotent group
(see \S 2.1 of \cite{SGA}). 
Base change for N\'eron models induces  
a homomorphism  
     $\iota:X_{v} \to X_{w}$
such that 
if $n$ is a positive integer not divisible by $p$, then the restriction of 
$\iota$ to the $n$-torsion $(X_{v})_n$ is injective
(see  Lemma 2 of \cite{SerreTate}; see also (3.1.1)
of \cite{SGA}). 
Here, $X_v$ and $X_w$ are viewed as commutative algebraic groups
over an algebraic closure of the residue field at $w$.
The map $\iota$ induces homomorphisms
$\T_{v} \to\T_{w}$ and $\B_{v} \to \B_{w}$
whose kernels are finite group schemes killed by 
appropriate powers of $p$.
The image of $\T_{v}$ (resp., $\B_v$) is a subtorus 
(resp., abelian subvariety) in $\T_{w}$
(resp., $\B_w$),
and we let $\T$ 
(resp., $\B$) denote the corresponding quotient.
Let $a$ and $t$ (respectively, $a_{v}$ and $t_{v}$)
denote the abelian and toric ranks of $X_{w}$ 
(respectively, $X_{v}$).
Note that $a$ and $t$ are independent of the valuation $w$ 
above $v$ at which $X$ has semistable reduction. 
We have $\rk(T_\ell(X_w))=2a+t$ and 
$\rk(T_\ell(X_v))=2a_v+t_v$.
By the functoriality of N\'eron models, $G_{v,X}$ acts on $X_{w}$
(see \S 4.2 of \cite{SGA}),
and therefore acts on  $\T_{w}$ and on $\B_w$.
One may easily check that this action is trivial on the 
image of $X_v \to X_w$.
It follows that $G_{v,X}$ acts on $\T$ and on $\B$.

Fix a polarization $\lambda$ on $X$. 
Let $W_{\ell}$ (respectively, $S_{\ell}$) 
denote the orthogonal complement of 
$V_{\ell}(X_{v})/V_{\ell}(\T_{v})$
(respectively, $T_{\ell}(X_{v})/T_{\ell}(\T_{v})$)
with respect to the 
pairing $e_{\lambda}$ on $V_{\ell}(X_{w})/V_{\ell}(\T_{w})$
(respectively, $T_{\ell}(X_{w})/T_{\ell}(\T_{w})$)
induced by $E_{\lambda}$. 
Then 
$W_{\ell}$ is a $G_{v,X}$-stable
$\Q_\ell$-vector space of dimension $2a-2a_{v}$ and
$S_{\ell}$ is a $G_{v,X}$-stable $\Z_\ell$-sublattice of 
rank $2a-2a_{v}$.

Recall that $\ell$ is always a prime not equal to $p$.

\begin{thm}
\label{hd}
\begin{enumerate}
\item[(i)] The form $e_\lambda: W_{\ell} \times W_{\ell} \to \Q_\ell$ is 
non-degenerate, alternating, and $G_{v,X}$-invariant.
The vector space $W_{\ell}$ and the lattice $S_{\ell}$ 
do not depend on the choice of 
polarization $\lambda$. 
The natural actions of $G_{v,X}$ on $\T$, 
$W_{\ell}$, and $e_\lambda$ induce an injection
$$G_{v,X} \hookrightarrow 
\Aut(\T) \times \Aut(W_{\ell},e_{\lambda}) \cong 
\GL_{t-t_{v}}(\Z) \times \Sp_{2(a-a_{v})}(\Q_{\ell})$$
such that 
the projection onto the first factor is independent of $\ell$,
and the characteristic polynomial of the projection onto the
second factor has integer coefficients independent of $\ell$.
\item[(ii)] 
If 
$\ell$ does not divide $\deg(\lambda)\#G_{v,X}$, then 
$e_\lambda: S_{\ell} \times S_{\ell} \to \Z_\ell$ is perfect and
the above injection takes values in 
$$\Aut(\T) \times \Aut(S_{\ell},e_{\lambda}) \cong 
\GL_{t-t_{v}}(\Z) \times \Sp_{2(a-a_{v})}(\Z_{\ell}).$$
\item[(iii)] 
Suppose $\ell$ does not divide $\#G_{v,X}$. Then 
for every $G_{v,X}$-stable $\Z_{\ell}$-lattice $\M$ in $W_{\ell}$, 
 there exists
 a perfect $\Z_{\ell}$-valued $G_{v,X}$-invariant alternating
pairing $e'$ on $\M$. I.e., there is an injection
$$G_{v,X} \hookrightarrow 
\Aut(\T) \times \Aut(\M,e') \cong 
\GL_{t-t_{v}}(\Z) \times \Sp_{2(a-a_{v})}(\Z_{\ell}).$$
\end{enumerate}
\end{thm}

\begin{proof}
By Proposition \ref{Wunipot}ii, $G_{v,X}$ acts faithfully on 
$T_{\ell}(X_{w})$.
The natural homomorphism
$$\varphi_{\ell} : G_{v,X} \to 
\Aut(T_{\ell}(\T_{w})) \times \Aut(T_{\ell}(X_{w})/T_{\ell}(\T_{w}))$$
is injective, since its kernel 
consists of unipotent operators of finite order
in characteristic zero.
The map $\varphi_{\ell}$ factors through 
$\Aut(\T) \times \Aut(T_{\ell}(X_{w})/T_{\ell}(\T_{w}))$,
since $G_{v,X}$ acts trivially on $\T_{v}$.
By Theorem \ref{serretateeqn}ii and Proposition \ref{orthogeqn}, 
$$T_{\ell}(\T_{v}) \subseteq 
T_{\ell}(X_{v}) \cap T_{\ell}(\T_{w}) =
T_{\ell}(X_{v}) \cap T_{\ell}(X_{w})^{\perp}
\subseteq 
T_{\ell}(X_{v}) \cap T_{\ell}(X_{v})^{\perp} =
T_{\ell}(\T_{v}).$$
Therefore,
$$T_{\ell}(\T_{v}) =
T_{\ell}(X_{v}) \cap T_{\ell}(\T_{w}) =
T_{\ell}(X_{v}) \cap T_{\ell}(X_{v})^{\perp}.$$
By Proposition \ref{orthogeqn}, 
$T_{\ell}(\T_{w}) = T_{\ell}(X_{w})^{\perp}$.
It follows that $e_\lambda$ 
is non-degenerate
on $T_\ell(X_w)/T_\ell(\T_w)$ and on
$T_\ell(X_v)/T_\ell(\T_v)$, and therefore on
$V_\ell(X_w)/V_\ell(\T_w)$ and on
$V_\ell(X_v)/V_\ell(\T_v)$. Further, $e_\lambda$ is
$G_{v,X}$-invariant.

Let $r = \#G_{v,X}$ and let
$u = \frac{1}{r}\sum_{g\in G_{v,X}}g \in \Q_\ell[G_{v,X}]$.
Then 
$$u(V_\ell(X_w)) = V_\ell(X_w)^{G_{v,X}} = V_\ell(X_v), \quad
u(V_\ell(\T_w)) = V_\ell(\T_w)^{G_{v,X}} = V_\ell(\T_v),$$
$$W_{\ell} = (1-u)(V_\ell(X_w)/V_\ell(\T_w)), \quad \text{ and }
S_{\ell}=W_{\ell} \cap (T_\ell(X_w)/T_\ell(\T_w)).$$
Therefore $W_{\ell}$ and $S_{\ell}$ are
independent of the choice of $\lambda$, and we have
a $G_{v,X}$-invariant $e_\lambda$-orthogonal splitting
$$V_{\ell}(X_{w})/V_{\ell}(\T_{w}) \cong 
V_{\ell}(X_{v})/V_{\ell}(\T_{v}) \oplus W_{\ell}.$$
Since $G_{v,X}$ acts trivially on $V_{\ell}(X_{v})/V_{\ell}(\T_{v})$, 
the map
$\varphi_{\ell}$ factors through 
$\Aut(\T) \times \Aut(W_{\ell})$.
Further, since $e_{\lambda}$ is non-degenerate, alternating, and
$G_{v,X}$-invariant, $\varphi_{\ell}$ factors through 
$\Aut(\T) \times \Aut(W_{\ell},e_{\lambda})
\cong \GL_{t-t_{v}}(\Z) \times \Sp_{2(a-a_{v})}(\Q_{\ell})$.
See p.~360 of \cite{SGA} for the integrality and
$\ell$-independence.

If $r$ is not divisible by $\ell$, then $u \in 
\Z_\ell[G_{v,X}]$, 
$$u(T_\ell(X_w)) = T_\ell(X_w)^{G_{v,X}} = T_\ell(X_v), \quad
u(T_\ell(\T_w)) = T_\ell(\T_w)^{G_{v,X}} = T_\ell(\T_v),$$
$$S_{\ell} = (1-u)(T_\ell(X_w)/T_\ell(\T_w)),$$  and we have
 a $G_{v,X}$-invariant $e_\lambda$-orthogonal splitting
$$T_{\ell}(X_{w})/T_{\ell}(\T_{w}) \cong 
T_{\ell}(X_{v})/T_{\ell}(\T_{v}) \oplus S_{\ell}.$$
To derive (iii), 
apply Proposition \ref{linalg}
and Remark \ref{perf}.
To derive (ii), 
suppose $\deg(\lambda)$ is not divisible by $\ell$. 
Then $E_{\lambda}$ is perfect.
Applying Proposition \ref{linalg1} with
$\M = T_\ell(X)$, $\N = T_\ell(X_w)$, $R = \Z_\ell$,
and $e = E_\lambda$, we deduce that $e_{\lambda}$ is
perfect, giving (ii).
\end{proof}

Next we give a variation of Theorem \ref{hd}, whose proof
and statement are of independent interest. 
Retain the notation from the proof of Theorem \ref{hd}.
Since $G_{v,X}$ is finite, there exists a 
$G_{v,X}$-invariant polarization $\delta$ on the abelian variety $\B$. 

\begin{thm}
\label{hd3}
If $\ell$ does not divide $\deg(\delta)$, then 
the action of $G_{v,X}$ on $\T$ and on $\B$ induces
injections
$$G_{v,X} 
\hookrightarrow \Aut(\T) \times \Aut(\B,\delta)
\hookrightarrow 
\GL_{t-t_{v}}(\Z) \times \Sp_{2(a-a_{v})}(\Z_{\ell})$$
where the projections
onto the first factors are independent of $\ell$,
and the characteristic polynomials of the projections onto the
second factors have integer coefficients independent of $\ell$.
In particular, the conclusion holds for all $\ell \ne p$,
whenever $a-a_v=0$ or $1$ 
(e.g., if $X$ or $\B$ is an elliptic curve).
\end{thm}

\begin{proof}
As $\I_{v}$-modules,  
$$T_{\ell}(\B_{w}) \cong T_{\ell}(X_{w})/T_{\ell}(\T_{w}), \quad
T_{\ell}(\B_{v}) \cong T_{\ell}(X_{v})/T_{\ell}(\T_{v}),$$
$$T_{\ell}(\B) \cong T_{\ell}(\B_{w})/T_{\ell}(\B_{v}), \quad 
\text{ and } \quad 
T_{\ell}(\T) \cong T_{\ell}(\T_{w})/T_{\ell}(\T_{v}),$$ 
and we have  
$\dim(\T) = t - t_{v}$ and $\dim(\B) = a - a_{v}$. 
By considering the pairing $E_\delta$ 
on $T_{\ell}(\B)$ induced by $\delta$, 
it follows that the image of $G_{v,X}$ in 
$\Aut(T_{\ell}(\B))$ lies in $\Aut(T_{\ell}(\B),E_\delta)$.
Since $\ell$ does not divide $\deg(\delta)$,
then $E_\delta$ is perfect.
The actions of $G_{v,X}$ on $\T$ and $\B$ therefore induce
a homomorphism 
$$\eta_{\ell} : G_{v,X} \to 
\Aut(\T) \times \Aut(T_\ell(\B),E_\delta) \cong 
\GL_{t-t_{v}}(\Z) \times 
\Sp_{2(a-a_{v})}(\Z_{\ell}).$$
By Proposition \ref{Wunipot}ii, $G_{v,X}$ acts faithfully on 
$T_{\ell}(X_{w})$.
The kernel of $\eta_{\ell}$ consists of elements 
   in     $G_{v,X} \subset \Aut(T_{\ell}(X_w))$
 which act as the identity on 
$T_{\ell}({\T}_w)$ and on $T_{\ell}(X_w)/T_{\ell}({\T}_w)$.
Therefore, these elements act as  unipotent operators on 
$T_{\ell}(X_w)$. Since they also have finite order, 
it follows that $\eta_{\ell}$ is injective.  
As before, see p.~360 of \cite{SGA} for  
integrality and $\ell$-independence.
When $\dim(\B)=0$ or $1$, we may suppose 
that $\deg(\delta)=1$.
\end{proof}

\begin{rem}
\label{finrem}
If $\ell \ne 2$, 
then reducing the second factor modulo $\ell$
in Theorem \ref{hd}iii or \ref{hd3} gives  
$$G_{v,X} \hookrightarrow 
\GL_{t-t_{v}}(\Z) \times \Sp_{2(a-a_{v})}(\F_{\ell}).$$
If either $p\#G_{v,X}$
or $p \deg(\delta)$ is odd, then
$$G_{v,X} \hookrightarrow 
\GL_{t-t_{v}}(\Z) \times \Sp_{2(a-a_{v})}(\Z/4\Z).$$
\end{rem}

\begin{rem}
It follows from Corollary \ref{hdcor} that the statement
in Theorem \ref{hd}iii
holds true whenever $p \ne \ell > 1+\max\{t-t_v,2(a-a_v)\}$
or $p \ne \ell > 2d+1$.
\end{rem}

\begin{rem}
Theorems \ref{hd} and \ref{hd3} and 
Remark \ref{finrem} have led us to the
question of when a finite interia group embedded in 
$\Sp_{2D}(\Q_{\ell})$ can also be embedded in 
$\Sp_{2D}(\Z_{\ell})$ or in $\Sp_{2D}(\F_{\ell})$
in such a way that the characteristic polynomials are
preserved.
Suppose $p$ is a prime number, and $G$ is a finite group with
a normal Sylow-$p$ subgroup $P$ such that $G/P$ is cyclic.
Suppose $\ell$ is a prime number, $\ell \ne p$, $V$ is a
$2D$-dimensional $\Q_{\ell}$-vector space, $e: V \times V \to \Q_{\ell}$
is a non-degenerate alternating bilinear form, and
$\varphi: G \hookrightarrow \Aut(V,e)$ is an injection.
Does there always exist a $G$-stable $\Z_{\ell}$-lattice $S$ in $V$
with a perfect $G$-invariant alternating $\Z_{\ell}$-valued pairing?
The answer is no. However, the answer is yes if $\ell > D+1$,
and this bound is sharp. 
Does there always exist an injection 
$\eta: G \hookrightarrow \Sp_{2D}(\F_{\ell})$
such that for all $g \in G$, the characteristic polynomial of
$\eta(g)$ is equal to the characteristic polynomial of
$\varphi(g)$ mod $\ell$?
The answer is no, but is yes if $\ell > 3$,
and this bound is sharp. 
Proofs and counterexamples (including counterexamples
with inertia groups of the
form $G_{v,X}$ for abelian surfaces $X$)
will appear in a later paper.
\end{rem}

\begin{rem}
Note that if $\#G_{v,X}$ is not divisible by $p$,
then $X$ acquires semistable reduction over
a tamely ramified extension;  
see \cite{Edixhoven} for a study of N\'eron models in 
this important case.
\end{rem}

\section{Bounds on the order of $G_{v,X}$}
\label{order}

One can use Theorem \ref{hd}, Remark \ref{cyclicrem}, and 
group theory to obtain more precise information about the
finite group $G_{v,X}$. In Corollary \ref{hdcor} below we
give one such result.
The next two results, along with Theorem \ref{hd}, 
will be used to prove Corollary \ref{hdcor}.

\begin{lem}
\label{2bd}
If $\ell$ is a prime number, $\ell \equiv 5 \pmod{8}$, 
$r$ and $m$ are positive integers, and 
$\Sp_{2m}(\F_{\ell})$ has an element of (exact) order $2^{r}$,
then $2^{r-1} \le 2m$.
\end{lem}

\begin{proof}
We may assume $r \ge 3$. 
Let $\zeta$ be a primitive $2^{r}$-th root of unity in 
${\bar \F_{\ell}}$. It is easy to check that the condition 
$\ell \equiv 5 \pmod{8}$ implies that 
$[\F_{\ell}(\zeta):\F_{\ell}] = 2^{r-2}$
and that $\zeta^{-1}$ is not a conjugate of $\zeta$.
It follows that in $\F_{\ell}[x]$ we have
$\Phi_{2^{r}} = fg$, where $\Phi_{2^{r}}$ is the $2^{r}$-th
cyclotomic polynomial, $f$ and $g$ are irreducible polynomials
in $\F_{\ell}[x]$ of degree $2^{r-2}$, and the roots of $g$ are
the inverses of the roots of $f$.
Let $\gamma$ denote an element of order $2^{r}$ in 
$\Sp_{2m}(\F_{\ell})$. If $\alpha$ is an eigenvalue of $\gamma$,
then so is $\alpha^{-1}$. It follows that the characteristic
polynomial of $\gamma$ is divisible by $\Phi_{2^{r}}$.
Therefore, $2^{r-1} \le 2m$.
\end{proof}

\begin{prop}
\label{qbd}
If $q$ is a prime number, $r$ and $m$ are positive integers,
and for all prime numbers $\ell$ in a set of density $1$
the group $\Sp_{2m}(\F_{\ell})$ contains an element of (exact) 
order $q^{r}$,
then $\varphi(q^{r}) \le 2m$.
\end{prop}

\begin{proof}
If $q = 2$ the result follows from Lemma \ref{2bd}.
Suppose $q$ is odd. Then $(\Z/q^{r}\Z)^{\times}$ is cyclic, and by
the Chebotarev density theorem there is a set of primes $\ell$
of positive density such that 
$\Gal(\F_{\ell}(\zeta_{q^{r}})/\F_{\ell}) \cong
(\Z/q^{r}\Z)^{\times} \cong
\Gal(\Q(\zeta_{q^{r}})/\Q)$.
If $q \ne \ell$, then
$\GL_{n}(\F_{\ell})$ contains an element of order $q^{r}$
if and only if $[\F_{\ell}(\zeta_{q^{r}}):\F_{\ell}] \le n$
(see \cite{Volvacev}).
Therefore, $\varphi(q^{r}) \le 2m$.
\end{proof}

An explicit description of all possible orders of elements of
general linear groups over arbitrary fields is given in 
\cite{Volvacev}.

Let $[\,\,\,]$ denote the greatest integer function,
let 
$s(n,q) = \sum_{j=0}^\infty \left[\frac{n}{q^j(q-1)}\right]$, 
and let 
$J(n) = \prod q^{s(n,q)}$, where $q$ runs over the prime numbers.
Note that the prime divisors of $J(n)$ are the primes 
$q \le n+1$.  
For example, $J(0) = 1$, $J(1) = 2$, and $J(2) = 24$.
For $n \ge 1$, Theorem 3.2 of \cite{JPAA} shows that
$J(n) < (4.462n)^{n}$ if $n$ is even and
$J(n) < \sqrt{2}(4.462n)^{n}$ if $n$ is odd.
The method of Minkowski and Serre (\cite{Mink} and pp.~119--121
of \cite{Serrearith}; see also Formula 3.1 of \cite{JPAA}) shows that,
for all $N \ge 3$, 
$J(2m)$ is equal to the greatest common divisor of the 
orders of the groups  $\Sp_{2m}(\F_{\ell})$, for primes $\ell \ge N$.
Further (\cite{SerreHarvard}), $J(n)$ is the least common multiple of 
the orders of the finite subgroups of $\GL_{n}(\Q)$ (or equivalently,
of $\GL_{n}(\Z)$).
While $J(n)$ is optimal from the point of view of divisibility,
there are sharper upper bounds on the orders of finite subgroups of
$\GL_{n}(\Q)$. 
The determination of the finite subgroups of maximum order for
general linear groups over $\Q$ and over cyclotomic fields is
given in \cite{Feit}.

Let
$$r_{p} = s(t-t_{v},p) + s(2(a-a_{v}),p), \qquad 
M_{v,X} = \max\{t-t_{v},2(a-a_{v})\},$$ and 
for all primes $q$ such that $p \ne q \le M_{v,X}+1$ let 
$$r_{q} = 1 + \left[\log_{q}\left(\frac{M_{v,X}}{q-1}\right)\right].$$
If $X$ has semistable reduction at $v$ let $N_{v,X} = 1$,
and otherwise let $N_{v,X} = \prod q^{r_{q}}$, where the
product runs over all prime numbers $q \le M_{v,X}+1$
(this might include $q=p$). 
Let $Q_{v,X}$ denote the largest prime divisor of 
$\#G_{v,X}$
(let $Q_{v,X} = 1$ if $\I_{v} = \II$).

\begin{cor}
\label{hdcor}
The order of $G_{v,X}$ divides $N_{v,X}$. In particular,
$Q_{v,X} \le M_{v,X} + 1 \le 2d + 1$, and 
$\#G_{v,X}$ divides 
$J(t-t_{v})J(2(a-a_{v}))$ and divides $J(2d)$. 
\end{cor}

\begin{proof}
By Theorem \ref{hd}ii,
 $\#G_{v,X}$ divides $J(t-t_{v})J(2(a-a_{v}))$,
and therefore $Q_{v,X} \le M_{v,X}+1$.
Note that $J(t-t_{v})J(2(a-a_{v}))$ divides $J(2d)$,
since 
$t - t_{v} + 2(a - a_{v}) \le t + 2a = 2d - t \le 2d$.
As noted in Remark \ref{cyclicrem}, the prime-to-$p$ part of 
$G_{v,X}$ is cyclic. 
Suppose $q$ is a prime divisor of $\#G_{v,X}$, and $q \ne p$.
Then $\GL_{n}(\Q)$ contains an element
of order $q^{r}$ if and only if $\varphi(q^{r}) \le n$, i.e.,
if and only if 
$r \le 1+\left[\log_{q}\left(\frac{n}{q-1}\right)\right]$.
The result now follows from Proposition \ref{qbd}.
\end{proof}

\section{Applications}
\label{applic}

Retain the notation of the previous sections. 
The next result follows immediately from 
Corollary \ref{hdcor} and Theorem \ref{main}iii,iv. 

\begin{cor}
\label{boundcor}
Suppose $L$ is a finite separable extension of $F$, and 
$w$ is the restriction of ${\bar v}$ to $L$. 
Suppose $X$ has semistable reduction at $w$ but not at $v$.
Then $[\I_{v}:\I_{w}]$ has a prime divisor $q$ such that 
$q \le Q_{v,X} \le M_{v,X} + 1 \le 2d + 1$.
\end{cor}

\begin{rem}
\label{divby}
Suppose $L$ is a finite separable extension of $F$, and
$w$ is the restriction of ${\bar v}$ to $L$.
Let $k_{w}$ and $k_{v}$ denote the residue fields and
let $e(w/v) = [w(L^\times):v(F^\times)]$ 
($=$ the ramification degree).
Then $[\I_{v}:\I_{w}] = e(w/v)[k_{w}:k_{v}]_{i}$, where
the subscript $i$ denotes the inseparable degree (see
Proposition 21 on p.~32 of \cite {Corps} for the case where
$L/F$ is Galois. In the non-Galois case, take a Galois extension
$L'$ of $F$ which contains $L$, and apply the result to
$L'/L$ and $L'/F$, to obtain the result for $L/F$).
Taking completions, then $[L_{w}:F_{v}] = e(w/v)[k_{w}:k_{v}]
= [\I_{v}:\I_{w}][k_{w}:k_{v}]_{s}$, where
the subscript $s$ denotes the separable degree.
Therefore, $[\I_{v}:\I_{w}]$ divides $[L_{w}:F_{v}]$.
\end{rem}

\begin{cor}
\label{prpow}
Suppose $L$ is a finite separable extension of $F$. 
Suppose in addition
that either $F$ is complete with respect to $v$,
or $L/F$ is Galois.
Suppose $X$ has semistable reduction at the restriction $w$  
of ${\bar v}$ to $L$, but does not have semistable reduction at $v$.
Then $[L:F]$ has a prime divisor $q$ such that $q \le Q_{v,X}$.
In particular, if $[L:F]$ is a power of a prime $q$, then 
$q \le Q_{v,X}$.
\end{cor}

\begin{proof}
Under our assumptions on $L/F$, its degree  
is divisible by $[I_{v}:\I_{w}]$. The result now follows from
Corollary \ref{boundcor}.
\end{proof}

Recall that there exists a finite Galois extension $L$ of $F$
such that $X$ has semistable reduction at the extensions 
of ${v}$ to $L$ (see Proposition 3.6 of \cite{SGA}).

\begin{cor}
\label{Gal}
Let $r = \#G_{v,X}$ and let $\zeta_{r}$ denote a primitive $r$-th 
root of unity.
Suppose that $r$ is not divisible by $p$.
Then there is a cyclic 
degree $r$ extension $L$ of $F(\zeta_{r})$ such
that $X$ acquires semistable reduction over every extension of $v$
to $L$.  
If either $F = F(\zeta_r)$ or $p > Q_{v,X}$, 
then there exists a finite Galois extension $L$ of $F$ 
of degree prime to $p$ such
that $X$ acquires semistable reduction over every extension of $v$
to $L$.  
\end{cor}

\begin{proof}
Let $L$ be the field obtained by adjoining 
an $r$-th root of a uniformizing parameter to $F(\zeta_{r})$.
Then $L/F$ is Galois,
$F(\zeta_{r})/F$ is unramified above $v$, $L/F(\zeta_{r})$
is totally ramified, and 
$\Gal(L/F(\zeta_{r})) \cong G_{v,X} \cong \Z/r\Z$.
Let $w$ be the restriction of ${\bar v}$ to $L$. 
By the construction of $L$, we have $\I_{w} \cong \II$. 
By Theorem \ref{main}iii, $X$ has semistable reduction at $w$.
Since $L/F$ is Galois, $X$ has semistable reduction at every
extension of $v$ to $L$ by Theorem \ref{galcrit}.
If $p > Q_{v,X}$, then  
$p$ does not divide $\varphi(r)$
by Corollary \ref{hdcor}, and therefore
$p$ does not divide $[F(\zeta_r):F]$. 
\end{proof}

By Corollary \ref{hdcor}, $Q_{v,X}$ can be replaced by $2d+1$ 
(or by $M_{v,X} + 1$) in Corollaries \ref{prpow} and
\ref{Gal}.


\begin{thebibliography}{99}

\bibitem{BLR} S.\ Bosch, W.\ L\"utkebohmert, M.\ Raynaud, N\'eron models,
Springer, Berlin-Heidelberg-New York, 1990.
\bibitem{Bourbaki} N. Bourbaki, Alg\`ebre, 
Chapitre~9, Hermann, Paris, 1959.
\bibitem{Edixhoven} B.\ Edixhoven, {\em N\'eron models and tame
ramification}, Comp.\ math.\ {\bf 81} (1992), 291--306.
\bibitem{Feit} W.\ Feit, {\em Orders of finite linear groups}, preprint.
\bibitem{Frohlich} A. Fr\"ohlich, {\em Local fields}, in Algebraic
Number Theory, J. W. S. Cassels and A. Fr\"ohlich, eds., Thompson
Book Company, Washington, 1967, pp. 1--41.
\bibitem{SGA} A.\ Grothendieck, {\em Mod\`eles de N\'eron et monodromie},
in Groupes de monodromie en g\'eometrie alg\'ebrique, SGA7 I 
(A.\ Gro\-then\-dieck, ed.), Lecture Notes in Math.\ {\bf 288}, Springer,
Berlin-Heidelberg-New York, 1972, pp.\ 313--523.
\bibitem{Lorenzini} D.\ Lorenzini, {\em
Groups of components of N\'eron models of
Jacobians}, Comp.\ math.\ 73 (1990), 145--160.
\bibitem{Mink} H.\ Minkowski, Gesammelte Abhandlungen, Bd.\ I, Leipzig, 
1911, pp.\ 212-218 ({\em Zur Theorie der positiven quadratischen Formen}, 
J.\ reine angew.\ Math.\ {\bf 101} (1887), 196--202).
\bibitem{Corps} J-P.\ Serre, Corps locaux, Hermann, Paris, 1968.
\bibitem{Serre72} J-P.\ Serre, {\em Propri\'et\'es galoisiennes des
points d'ordre fini des courbes elliptiques}, Invent.\ math.\ {\bf 15}
(1972), 259--331.
\bibitem{serrereps} J-P.\ Serre, {\em Repr\'esentations
lin\'eaires des groupes finis}, Third edition, Hermann, Paris,
1978.
\bibitem{Serrearith} J-P.\ Serre, {\em Arithmetic Groups},
in Homological Group Theory (C.\ T.\ C.\ Wall, ed.),
LMS Lecture Note Series {\bf 36}, Camb.\ Univ.\ Press,
Cambridge, 1979, pp.\ 105--136.
\bibitem{Motives} J-P.\ Serre, {\em Propri\'et\'es conjecturales 
des groupes de Galois motiviques et des repr\'esentations $\ell$-adiques}, 
in Motives 
(U.\ Jannsen, S.\ Kleiman, J-P.\ Serre, eds.), Proc.\ Symp.\ Pure 
Math.\ {\bf 55} (1994), Part 1, pp.~377--400.
\bibitem{SerreHarvard} J-P.\ Serre, {\em Finite subgroups of Lie 
groups}, Lectures at Harvard, Autumn 1994, Lecture~2. 
\bibitem{SerreTate} J-P.\ Serre, J.\ Tate, {\em Good reduction of abelian 
varieties}, Ann.\ of Math.\ {\bf 88} (1968), 492--517. 
\bibitem{JPAA} A.\ Silverberg,  {\em Fields of definition for homomorphisms 
of abelian varieties}, J.\ Pure and Applied Algebra {\bf 77} (1992),
253--262.
\bibitem{semistab} A.\ Silverberg, Yu.\ G.\ Zarhin, {\em Semistable reduction 
and torsion subgroups of abelian varieties}, 
Ann.\ Inst.\ Fourier {\bf 45}, no.~2 (1995), 403--420.
\bibitem{Compositio} A.\ Silverberg, Yu.\ G.\ Zarhin, {\em 
Connectedness results for $\ell$-adic 
representations associated to abelian varieties}, 
Comp.\ math.\ {\bf 97} (1995), 273--284.
\bibitem{smalldeg} A.\ Silverberg, Yu.\ G.\ Zarhin, 
{\em Semistable reduction of abelian varieties over extensions
of small degree}, to appear in J.\ Pure and Applied Algebra.
\bibitem{Volvacev} R.\ T.\ Vol'va\v{c}ev, {\em On the order of an element
of a matrix group} (Russian), Vesci Akad.\ Navuk BSSR Ser.\ Fiz.-Mat.\
Navuk no.\ 2 (1965), 11--16 ({\bf MR~34}~\#{5943}).
\end{thebibliography}
\end{document}